# An Analysis of Bennett's Pebble Game

E. Knill, knill@lanl.gov*

May 1995

**Abstract**

Bennett's pebble game [1, 2] was introduced to obtain better time/space tradeoffs in the simulation of standard Turing machines by reversible ones. So far only upper bounds for the tradeoff based on the pebble game have been published. Here we give a recursion for the time optimal solution of the pebble game given a space bound. We analyze the recursion to obtain an explicit asymptotic expression for the best time-space product.

## 1 Introduction

In Bennett's pebble game, we have a board of $n$ squares labeled $s_1, \ldots, s_n$. Each square may be empty or contain one pebble. A move in the game consists of adding or removing a pebble. A pebble can be added to or removed from square $s_i$ if $i = 1$ or if there is a pebble on square $s_{i-1}$. Initially, the board is empty. The desired final configuration has exactly one pebble on the board on square $s_n$. For the motivation and application to reversible computing, we refer the reader to [1, 2].

The pebble game is trivial without resource constraints. At the $i$'th step for $1 \leq i \leq n$ we can place a pebble on square $s_i$. At the $n + i$'th step for $1 \leq i \leq n - 1$ we can remove the pebble on square $s_{n-i}$. This solves the game in $2n - 1$ steps using $n$ pebbles. The time resource $T$ is the number of steps to solve the game, while the space resource $S$ is the number of pebbles used, that is the maximum number of pebbles on the board at any time. For any solution $\mathcal{S}$ of the pebble game, let $T(\mathcal{S})$ be the number of steps and $S(\mathcal{S})$ the number of pebbles used. Our goal is to determine

$$F(n, S) = \min(T(\mathcal{S}) \mid S(\mathcal{S}) \leq S).$$

Levine and Sherman [2] expressed Bennett's [1] upper bound on $F(n, S)$ by a recursion which they were able to solve. We will obtain an exact recursive expression for $F(n, S)$ which can be evaluated explicitly for any desired $n$ and $S$. An asymptotic expression for the minimum time space product is obtained. A consequence of this analysis is that an upper bound for the minimum time-space product of quantum factoring (see [3]) is $2^{2\sqrt{n}(1+o(1))} n^3$, where $n$ is the number of bits in the number to be factored. We defer discussion of this observation to another note.

## 2 Analysis of the pebble game

### 2.1 Solutions of the pebble game

To obtain a recursive expression for $F(n, S)$ it is useful to model solutions $\mathcal{S}$ by representing the residence of each pebble on a square as an interval. Formally, we view $\mathcal{S}$ as a $0 - 1$ matrix, also denoted by $\mathcal{S}$, of dimensions $n \times [-\infty, +\infty]$. $\mathcal{S}_{i,j} = 1$ means that there is a pebble on square $s_i$ after step $j$. For $\mathcal{S}$ to be a valid solution, we need

(i) Initial configuration: There exists $T_b$ such that for $t \leq T_b$, $\mathcal{S}_{i,t} = 0$.

(ii) Final configuration: There exists $T_e$, such that for $t \geq T_e$, $\mathcal{S}_{i,t} = \delta_{n,i}$, with $\delta_{n,i}$ the usual delta function.

(iii) One move at a time: If $\mathcal{S}_{i,j} \neq \mathcal{S}_{i,j+1}$ then for all $i' \neq i$, $\mathcal{S}_{i',j} = \mathcal{S}_{i',j+1}$.

(iv) Add constraint: If $\mathcal{S}_{i,j} = 0$ and $\mathcal{S}_{i,j+1} = 1$, then either $i = 1$ or $\mathcal{S}_{i-1,j} = 1$.

(v) Remove constraint: If $\mathcal{S}_{i,j} = 0$ and $\mathcal{S}_{i,j-1} = 1$, then either $i = 1$ or $\mathcal{S}_{i-1,j-1} = 1$.

The space requirement of $\mathcal{S}$ is given by the maximum number of 1's in any column of $\mathcal{S}$, while the time requirement is the number of inequalities between adjacent elements in the rows.

Consider square $s_i$. The $i$'th row of $\mathcal{S}$ can be viewed as a collection of intervals of 1's, representing the residence of a pebble on $s_i$. Thus a *(residence) interval* $[k,l]$ of $s_i$ satisfies (1) $k \leq l$, (2) $\mathcal{S}_{i,j} = 1$ for $k \leq j \leq l$, (3) $\mathcal{S}_{i,k-1} = 0$ and (4) $\mathcal{S}_{i,l+1} = 0$.

Let $j, j' \geq 0$. If $I = [k,l]$ is an interval of $s_i$ and $J = [k+j, l-j']$ is an interval of $s_{i-1}$, then $J$ does not contribute anything to the solution. A better solution in terms of time and space is obtained by removing $J$ (by setting the corresponding entries of $\mathcal{S}$ to 0). This allows us to add another constraint:

(vi) If $I$ is an interval of $s_i$ and $J$ is an interval of $s_{i-1}$, then $J \not\subseteq I$.

In depicting solutions it is convenient to represent the intervals as horizontal lines at a vertical offset, where no two endpoints have the same horizontal offset. See Figure 2.1 for two examples.

## 2.2 A recursion for $F(n, S)$

**Theorem 2.1** *$F(n, S)$ is given by*

$$F(n, S) = \begin{cases} 1 & \text{for } n = 1 \text{ and } S \geq 1, \\ \infty & \text{for } n \geq 2 \text{ and } S = 1, \\ \min(F(m, S) + F(m, S-1) + F(n-m, S-1) \mid 1 \leq m < n) & \text{for } n \geq 2 \text{ and } S \geq 2. \end{cases}$$

**Proof.** For $n = 1$, the solution is to place a pebble on $s_1$ in one step. For $n > 1$, at least 2 pebbles are required.

The recursive expression for $F(n, S)$ is obtained by decomposing the solution matrix. The basic idea is shown in Figure 2.1. We begin by proving that $F(n, S) \geq \min(F(m, S) + F(m, S-1) + F(n-m, S-1) \mid 1 \leq m < n)$. Let $n > 1$ and $S > 1$. Let $\mathcal{S}$ be a time optimal solution of the pebble game with $S(\mathcal{S}) \leq S$. Let $t_b$ be the least $t$ such that $\mathcal{S}_{1,t} = 1$. We can assume that for $t \geq t_b$, the $t$'th column of $\mathcal{S}$ always has at least one 1. If not we can obtain a more time efficient solution by setting all the entries preceding that column to 0. Let $t_e$ be the largest $t$ such that $\mathcal{S}_{1,t} = 1$. Due to the move constraints, for $t > t_e$, $\mathcal{S}_{i,t} = \delta_{n,i}$. Define $m(\mathcal{S}, t) = \min(i \mid \mathcal{S}_{i,t} = 1)$ and let $t_m$ be the maximum $t$ with $t_b \leq t \leq t_e$ for which $m(\mathcal{S}, t)$ achieves its maximum in this range. Let $m = m(\mathcal{S}, t_m)$. Observe that $m < n$, else we can obtain a better solution by setting $\mathcal{S}_{i,t} = 0$ for $i < m$ and $t > t_m$. We can extract three solutions of smaller problems from $\mathcal{S}$.

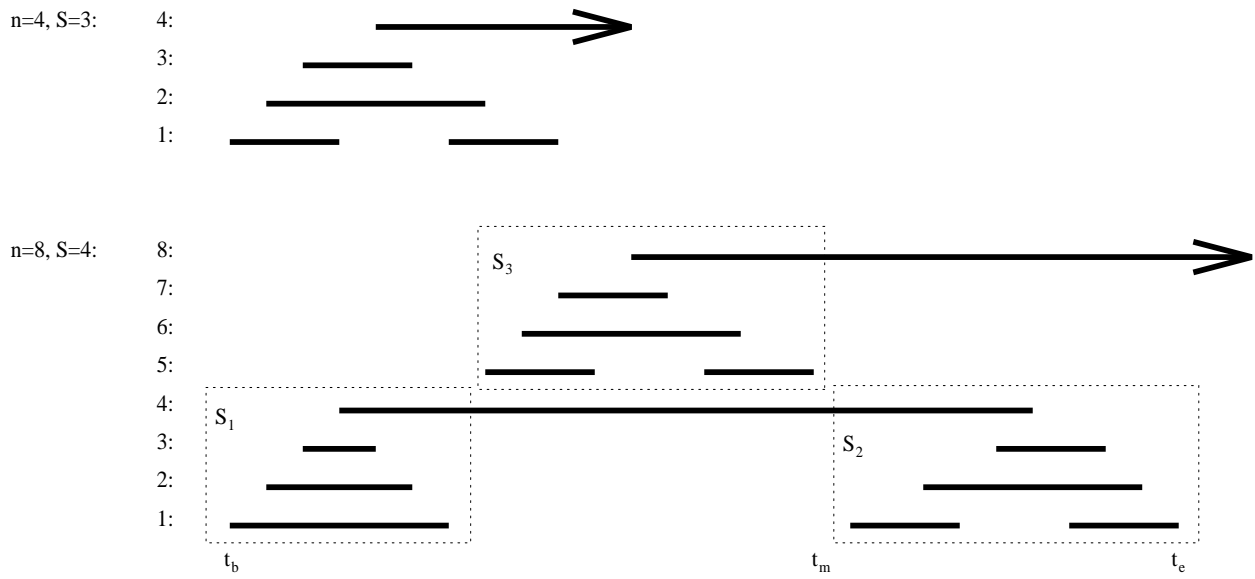

Figure 1. Solutions as residence intervals.

Let $\mathcal{S}'_1$ be the submatrix of $\mathcal{S}$ consisting of the first $m$ rows between columns $t_b$ and $t_m$. Let $\mathcal{S}_1$ be the extension of $\mathcal{S}'_1$ to a solution of the pebble game obtained by letting $\mathcal{S}_{1;i,t} = 0$ for $t < t_b$ and $\mathcal{S}_{1;i,t} = \delta_{m,i}$ for $t > t_m$. Then $\mathcal{S}_1$ is a solution of the pebble game with $m$ squares using at most $S$ pebbles.

Let $\mathcal{S}'_2$ be the the submatrix of $\mathcal{S}$ consisting of the first $m$ rows between columns $t_m$ and $t_e$. Let $\mathcal{S}_2$ be $\mathcal{S}'_2$ with the columns in reverse order and extended to a solution as above for $\mathcal{S}'_1$. Then $\mathcal{S}_2$ is a solution of the pebble game with $m$ squares. We claim that it uses at most $S - 1$ pebbles. To see this we show that for $t \geq t_m$, the maximum $i(t)$ for which $\mathcal{S}_{i,t} = 1$ satisfies $i > m$. Suppose that this is not the case. Let $t_1$ and $t_2$ be the right endpoint and the left endpoint of the two intervals of $s_1$ surrounding $t_m$. Because of the move constraints, $\mathcal{S}_{i,t} = 0$ for $t_1 < t < t_2$ and $i < m$. This implies that the columns between $t_1$ and $t_2$ have only one 1. Since $S \geq 2$, we can set $\mathcal{S}_{1,t}$ to 1 for $t_1 < t < t_2$. This reduces the number of steps by 2, contradicting optimality of $\mathcal{S}$. This argument shows that $i(t_m) > m$. If for $t' > t_m$, $i(t') \leq m$, then we can obtain a more efficient solution by setting $\mathcal{S}_{i,t} = 0$ for $i > m$ and $t \leq t'$.

Let $\mathcal{S}_3$ be the matrix consisting of the rows with index $i > m$. This matrix is a solution of the pebble game with $n - m$ squares. Since this solution achieves the final configuration before $t_e$, and by the choice of $m$, $\mathcal{S}_3$ can use at most $S - 1$ pebbles.

Since $T(\mathcal{S}) = T(\mathcal{S}_1) + T(\mathcal{S}_2) + T(\mathcal{S}_3)$, the desired inequality follows.

To show that $F(n, S) \leq \min(F(m, S) + F(m, S - 1) + F(n - m, S - 1) \mid 1 \leq m < n)$, we use the decomposition for the lower bound and show that it can be used constructively for an upper bound. Given $\mathcal{S}_1$, $\mathcal{S}_2$ and $\mathcal{S}_3$ with $S(\mathcal{S}_1) \leq S$, $S(\mathcal{S}_2) \leq S - 1$ and $S(\mathcal{S}_3) \leq S - 1$, we construct a solution for the pebble game with $n$ squares by composing the solutions $\mathcal{S}_i$. Let $t_{b,i}$ and $t_{e,i}$ be defined for $\mathcal{S}_i$ as we defined $t_b$ and $t_e$ for $\mathcal{S}$ above (if $m = 1$ or $n - m = 1$, set the corresponding $t_e = t_b + 1$). By shifting the solutions, we can assume $t_{e,1} + 2 = t_{b,3}$ and $t_{e,2} = t_{e,3} + 2$. Define

$$\mathcal{S}_{i,t} = \begin{cases} \mathcal{S}_{1;i,t} & \text{for } i \leq m \text{ and } t \leq t_{e,1} + 1, \\ \mathcal{S}_{3;i-m,t} & \text{for } m < i \leq n, \\ \mathcal{S}_{2;i,2t_{e,2}-t} & \text{for } i \leq m \text{ and } t_{e,3} + 1 < t, \\ 1 & \text{for } i = m \text{ and } t_{b,3} \leq t \leq t_{e,3} + 1, \\ 0 & \text{otherwise.} \end{cases}$$

Then $\mathcal{S}$ is a solution with $T(\mathcal{S}) = \sum_i T(\mathcal{S}_i)$ satisfying the space bound. ■

The recursion of Theorem 2.1 can be used to compute tables of $F(n, S)$.

The appendix contains some Mathematica expressions which compute $F(n, S)$ dynamically. Tables 1 and 2 contain $F(n, S)$ for small values of $n$ and $S$. The proof of the theorem can be used to easily obtain the optimal strategy for each $n$ and $S$.

## 2.3 Properties of $F(n, S)$

Some immediate observations are as follows:

**Observation 2.2** *(i) $F(n, S) = 2n - 1$ for $S \geq n$.*

*(ii) For $S \leq S'$, $F(n, S) \geq F(n, S')$.*

*(iii) For $n \leq n'$, $F(n, S) \leq F(n', S)$.*

We next establish the maximum $n$ for which $F(n, S) \leq \infty$.

**Theorem 2.3** $F(n, S) < \infty$ *iff* $n \leq 2^{S-1}$.

**Proof.** The statement holds for $S = 1$. We proceed by induction on $S$ and $n$. If $n \leq 2^{S-1}$ and we let $m = \lfloor n/2 \rfloor$, then $m \leq 2^{S-2}$ and $n - m \leq 2^{S-2}$. It follows that $F(m, S - 1)$, $F(n - m, S - 1)$ and $F(m, S)$ are finite by the induction hypothesis. The recursive expression for $F(n, S)$ now implies that $F(n, S)$ is finite.

If $n > 2^{S-1}$ then either $m > 2^{S-2}$ or $n - m > 2^{S-2}$. Hence every term in the minimum of the recursion for $F(n, S)$ is infinite. ∎

To obtain explicit bounds on $F(n, S)$ it is convenient to consider the function

$$\Delta(n, S) = \begin{cases} 0 & \text{if } n \leq 0, \\ F(n + 1, S) - F(n, S) & \text{if } n > 0 \text{ and } F(n + 1, S) \leq \infty, \\ \infty & \text{otherwise.} \end{cases}$$

For $n > 1$ and $S > 1$, let $m(n, S)$ be the least minimizer of the recursive expression for $F(n, S)$. Define $m(0, S) = m(1, S) = 0$. Note that $m(n, S) \leq \lfloor n/2 \rfloor$ by the monotonicity properties of $F(n, S)$.

**Lemma 2.4** . *For $m = m(n, S)$, $n \geq 1$ and $S \geq 2$,*

$$\Delta(n - m, S - 1) > \Delta(m - 1, S) + \Delta(m - 1, S - 1),$$
$$\Delta(n - m - 1, S - 1) \leq \Delta(m, S) + \Delta(m, S - 1).$$

| n \ S | 1 | 2 | 3 | 4 | 5 | 6 | 7 | 8 | 9 | 10 | 11 | 12 | 13 | 14 | 15 | 16 | 17 | 18 | 19 | 20 |
|---|---|---|---|---|---|---|---|---|---|---|---|---|---|---|---|---|---|---|---|---|
| 1 | 1 | 1 | 1 | 1 | 1 | 1 | 1 | 1 | 1 | 1 | 1 | 1 | 1 | 1 | 1 | 1 | 1 | 1 | 1 | 1 |
| 2 | ∞ | 3 | 3 | 3 | 3 | 3 | 3 | 3 | 3 | 3 | 3 | 3 | 3 | 3 | 3 | 3 | 3 | 3 | 3 | 3 |
| 3 | ∞ | ∞ | 5 | 5 | 5 | 5 | 5 | 5 | 5 | 5 | 5 | 5 | 5 | 5 | 5 | 5 | 5 | 5 | 5 | 5 |
| 4 | ∞ | ∞ | 9 | 7 | 7 | 7 | 7 | 7 | 7 | 7 | 7 | 7 | 7 | 7 | 7 | 7 | 7 | 7 | 7 | 7 |
| 5 | ∞ | ∞ | ∞ | 11 | 9 | 9 | 9 | 9 | 9 | 9 | 9 | 9 | 9 | 9 | 9 | 9 | 9 | 9 | 9 | 9 |
| 6 | ∞ | ∞ | ∞ | 15 | 13 | 11 | 11 | 11 | 11 | 11 | 11 | 11 | 11 | 11 | 11 | 11 | 11 | 11 | 11 | 11 |
| 7 | ∞ | ∞ | ∞ | 19 | 17 | 15 | 13 | 13 | 13 | 13 | 13 | 13 | 13 | 13 | 13 | 13 | 13 | 13 | 13 | 13 |
| 8 | ∞ | ∞ | ∞ | 25 | 21 | 19 | 17 | 15 | 15 | 15 | 15 | 15 | 15 | 15 | 15 | 15 | 15 | 15 | 15 | 15 |
| 9 | ∞ | ∞ | ∞ | ∞ | 25 | 23 | 21 | 19 | 17 | 17 | 17 | 17 | 17 | 17 | 17 | 17 | 17 | 17 | 17 | 17 |
| 10 | ∞ | ∞ | ∞ | ∞ | 29 | 27 | 25 | 23 | 21 | 19 | 19 | 19 | 19 | 19 | 19 | 19 | 19 | 19 | 19 | 19 |
| 11 | ∞ | ∞ | ∞ | ∞ | 33 | 31 | 29 | 27 | 25 | 23 | 21 | 21 | 21 | 21 | 21 | 21 | 21 | 21 | 21 | 21 |
| 12 | ∞ | ∞ | ∞ | ∞ | 39 | 35 | 33 | 31 | 29 | 27 | 25 | 23 | 23 | 23 | 23 | 23 | 23 | 23 | 23 | 23 |
| 13 | ∞ | ∞ | ∞ | ∞ | 45 | 39 | 37 | 35 | 33 | 31 | 29 | 27 | 25 | 25 | 25 | 25 | 25 | 25 | 25 | 25 |
| 14 | ∞ | ∞ | ∞ | ∞ | 53 | 43 | 41 | 39 | 37 | 35 | 33 | 31 | 29 | 27 | 27 | 27 | 27 | 27 | 27 | 27 |
| 15 | ∞ | ∞ | ∞ | ∞ | 61 | 47 | 45 | 43 | 41 | 39 | 37 | 35 | 33 | 31 | 29 | 29 | 29 | 29 | 29 | 29 |
| 16 | ∞ | ∞ | ∞ | ∞ | 71 | 51 | 49 | 47 | 45 | 43 | 41 | 39 | 37 | 35 | 33 | 31 | 31 | 31 | 31 | 31 |
| 17 | ∞ | ∞ | ∞ | ∞ | ∞ | 57 | 53 | 51 | 49 | 47 | 45 | 43 | 41 | 39 | 37 | 35 | 33 | 33 | 33 | 33 |
| 18 | ∞ | ∞ | ∞ | ∞ | ∞ | 63 | 57 | 55 | 53 | 51 | 49 | 47 | 45 | 43 | 41 | 39 | 37 | 35 | 35 | 35 |
| 19 | ∞ | ∞ | ∞ | ∞ | ∞ | 69 | 61 | 59 | 57 | 55 | 53 | 51 | 49 | 47 | 45 | 43 | 41 | 39 | 37 | 37 |
| 20 | ∞ | ∞ | ∞ | ∞ | ∞ | 77 | 65 | 63 | 61 | 59 | 57 | 55 | 53 | 51 | 49 | 47 | 45 | 43 | 41 | 39 |
| 21 | ∞ | ∞ | ∞ | ∞ | ∞ | 85 | 69 | 67 | 65 | 63 | 61 | 59 | 57 | 55 | 53 | 51 | 49 | 47 | 45 | 43 |
| 22 | ∞ | ∞ | ∞ | ∞ | ∞ | 93 | 73 | 71 | 69 | 67 | 65 | 63 | 61 | 59 | 57 | 55 | 53 | 51 | 49 | 47 |
| 23 | ∞ | ∞ | ∞ | ∞ | ∞ | 101 | 79 | 75 | 73 | 71 | 69 | 67 | 65 | 63 | 61 | 59 | 57 | 55 | 53 | 51 |
| 24 | ∞ | ∞ | ∞ | ∞ | ∞ | 109 | 85 | 79 | 77 | 75 | 73 | 71 | 69 | 67 | 65 | 63 | 61 | 59 | 57 | 55 |
| 25 | ∞ | ∞ | ∞ | ∞ | ∞ | 117 | 91 | 83 | 81 | 79 | 77 | 75 | 73 | 71 | 69 | 67 | 65 | 63 | 61 | 59 |
| 26 | ∞ | ∞ | ∞ | ∞ | ∞ | 125 | 97 | 87 | 85 | 83 | 81 | 79 | 77 | 75 | 73 | 71 | 69 | 67 | 65 | 63 |
| 27 | ∞ | ∞ | ∞ | ∞ | ∞ | 135 | 105 | 91 | 89 | 87 | 85 | 83 | 81 | 79 | 77 | 75 | 73 | 71 | 69 | 67 |
| 28 | ∞ | ∞ | ∞ | ∞ | ∞ | 145 | 113 | 95 | 93 | 91 | 89 | 87 | 85 | 83 | 81 | 79 | 77 | 75 | 73 | 71 |
| 29 | ∞ | ∞ | ∞ | ∞ | ∞ | 155 | 121 | 99 | 97 | 95 | 93 | 91 | 89 | 87 | 85 | 83 | 81 | 79 | 77 | 75 |
| 30 | ∞ | ∞ | ∞ | ∞ | ∞ | 167 | 129 | 105 | 101 | 99 | 97 | 95 | 93 | 91 | 89 | 87 | 85 | 83 | 81 | 79 |
| 31 | ∞ | ∞ | ∞ | ∞ | ∞ | 179 | 137 | 111 | 105 | 103 | 101 | 99 | 97 | 95 | 93 | 91 | 89 | 87 | 85 | 83 |
| 32 | ∞ | ∞ | ∞ | ∞ | ∞ | 193 | 145 | 117 | 109 | 107 | 105 | 103 | 101 | 99 | 97 | 95 | 93 | 91 | 89 | 87 |
| 33 | ∞ | ∞ | ∞ | ∞ | ∞ | ∞ | 153 | 123 | 113 | 111 | 109 | 107 | 105 | 103 | 101 | 99 | 97 | 95 | 93 | 91 |
| 34 | ∞ | ∞ | ∞ | ∞ | ∞ | ∞ | 161 | 129 | 117 | 115 | 113 | 111 | 109 | 107 | 105 | 103 | 101 | 99 | 97 | 95 |
| 35 | ∞ | ∞ | ∞ | ∞ | ∞ | ∞ | 169 | 137 | 121 | 119 | 117 | 115 | 113 | 111 | 109 | 107 | 105 | 103 | 101 | 99 |
| 36 | ∞ | ∞ | ∞ | ∞ | ∞ | ∞ | 177 | 145 | 125 | 123 | 121 | 119 | 117 | 115 | 113 | 111 | 109 | 107 | 105 | 103 |
| 37 | ∞ | ∞ | ∞ | ∞ | ∞ | ∞ | 185 | 153 | 129 | 127 | 125 | 123 | 121 | 119 | 117 | 115 | 113 | 111 | 109 | 107 |
| 38 | ∞ | ∞ | ∞ | ∞ | ∞ | ∞ | 193 | 161 | 135 | 131 | 129 | 127 | 125 | 123 | 121 | 119 | 117 | 115 | 113 | 111 |
| 39 | ∞ | ∞ | ∞ | ∞ | ∞ | ∞ | 201 | 169 | 141 | 135 | 133 | 131 | 129 | 127 | 125 | 123 | 121 | 119 | 117 | 115 |
| 40 | ∞ | ∞ | ∞ | ∞ | ∞ | ∞ | 209 | 177 | 147 | 139 | 137 | 135 | 133 | 131 | 129 | 127 | 125 | 123 | 121 | 119 |
| 41 | ∞ | ∞ | ∞ | ∞ | ∞ | ∞ | 217 | 185 | 153 | 143 | 141 | 139 | 137 | 135 | 133 | 131 | 129 | 127 | 125 | 123 |
| 42 | ∞ | ∞ | ∞ | ∞ | ∞ | ∞ | 225 | 193 | 159 | 147 | 145 | 143 | 141 | 139 | 137 | 135 | 133 | 131 | 129 | 127 |
| 43 | ∞ | ∞ | ∞ | ∞ | ∞ | ∞ | 235 | 201 | 165 | 151 | 149 | 147 | 145 | 143 | 141 | 139 | 137 | 135 | 133 | 131 |
| 44 | ∞ | ∞ | ∞ | ∞ | ∞ | ∞ | 245 | 209 | 173 | 155 | 153 | 151 | 149 | 147 | 145 | 143 | 141 | 139 | 137 | 135 |
| 45 | ∞ | ∞ | ∞ | ∞ | ∞ | ∞ | 255 | 217 | 181 | 159 | 157 | 155 | 153 | 151 | 149 | 147 | 145 | 143 | 141 | 139 |
| 46 | ∞ | ∞ | ∞ | ∞ | ∞ | ∞ | 265 | 225 | 189 | 163 | 161 | 159 | 157 | 155 | 153 | 151 | 149 | 147 | 145 | 143 |
| 47 | ∞ | ∞ | ∞ | ∞ | ∞ | ∞ | 275 | 233 | 197 | 169 | 165 | 163 | 161 | 159 | 157 | 155 | 153 | 151 | 149 | 147 |
| 48 | ∞ | ∞ | ∞ | ∞ | ∞ | ∞ | 285 | 241 | 205 | 175 | 169 | 167 | 165 | 163 | 161 | 159 | 157 | 155 | 153 | 151 |
| 49 | ∞ | ∞ | ∞ | ∞ | ∞ | ∞ | 297 | 249 | 213 | 181 | 173 | 171 | 169 | 167 | 165 | 163 | 161 | 159 | 157 | 155 |
| 50 | ∞ | ∞ | ∞ | ∞ | ∞ | ∞ | 309 | 257 | 221 | 187 | 177 | 175 | 173 | 171 | 169 | 167 | 165 | 163 | 161 | 159 |

Table 1. Table of $F(n, S)$ for $1 \leq n \leq 50$ and $1 \leq S \leq 20$.

| n \ S | 1 | 2 | 3 | 4 | 5 | 6 | 7 | 8 | 9 | 10 | 11 | 12 | 13 | 14 | 15 | 16 | 17 | 18 | 19 | 20 |
|---|---|---|---|---|---|---|---|---|---|---|---|---|---|---|---|---|---|---|---|---|
| 51 | ∞ | ∞ | ∞ | ∞ | ∞ | ∞ | 321 | 265 | 229 | 193 | 181 | 179 | 177 | 175 | 173 | 171 | 169 | 167 | 165 | 163 |
| 52 | ∞ | ∞ | ∞ | ∞ | ∞ | ∞ | 333 | 273 | 237 | 199 | 185 | 183 | 181 | 179 | 177 | 175 | 173 | 171 | 169 | 167 |
| 53 | ∞ | ∞ | ∞ | ∞ | ∞ | ∞ | 345 | 281 | 245 | 205 | 189 | 187 | 185 | 183 | 181 | 179 | 177 | 175 | 173 | 171 |
| 54 | ∞ | ∞ | ∞ | ∞ | ∞ | ∞ | 359 | 289 | 253 | 213 | 193 | 191 | 189 | 187 | 185 | 183 | 181 | 179 | 177 | 175 |
| 55 | ∞ | ∞ | ∞ | ∞ | ∞ | ∞ | 373 | 297 | 261 | 221 | 197 | 195 | 193 | 191 | 189 | 187 | 185 | 183 | 181 | 179 |
| 56 | ∞ | ∞ | ∞ | ∞ | ∞ | ∞ | 387 | 305 | 269 | 229 | 201 | 199 | 197 | 195 | 193 | 191 | 189 | 187 | 185 | 183 |
| 57 | ∞ | ∞ | ∞ | ∞ | ∞ | ∞ | 401 | 313 | 277 | 237 | 207 | 203 | 201 | 199 | 197 | 195 | 193 | 191 | 189 | 187 |
| 58 | ∞ | ∞ | ∞ | ∞ | ∞ | ∞ | 415 | 321 | 285 | 245 | 213 | 207 | 205 | 203 | 201 | 199 | 197 | 195 | 193 | 191 |
| 59 | ∞ | ∞ | ∞ | ∞ | ∞ | ∞ | 433 | 329 | 293 | 253 | 219 | 211 | 209 | 207 | 205 | 203 | 201 | 199 | 197 | 195 |
| 60 | ∞ | ∞ | ∞ | ∞ | ∞ | ∞ | 451 | 337 | 301 | 261 | 225 | 215 | 213 | 211 | 209 | 207 | 205 | 203 | 201 | 199 |
| 61 | ∞ | ∞ | ∞ | ∞ | ∞ | ∞ | 469 | 345 | 309 | 269 | 231 | 219 | 217 | 215 | 213 | 211 | 209 | 207 | 205 | 203 |
| 62 | ∞ | ∞ | ∞ | ∞ | ∞ | ∞ | 489 | 353 | 317 | 277 | 237 | 223 | 221 | 219 | 217 | 215 | 213 | 211 | 209 | 207 |
| 63 | ∞ | ∞ | ∞ | ∞ | ∞ | ∞ | 509 | 361 | 325 | 285 | 243 | 227 | 225 | 223 | 221 | 219 | 217 | 215 | 213 | 211 |
| 64 | ∞ | ∞ | ∞ | ∞ | ∞ | ∞ | 531 | 369 | 333 | 293 | 249 | 231 | 229 | 227 | 225 | 223 | 221 | 219 | 217 | 215 |
| 65 | ∞ | ∞ | ∞ | ∞ | ∞ | ∞ | ∞ | 379 | 341 | 301 | 257 | 235 | 233 | 231 | 229 | 227 | 225 | 223 | 221 | 219 |
| 66 | ∞ | ∞ | ∞ | ∞ | ∞ | ∞ | ∞ | 389 | 349 | 309 | 265 | 239 | 237 | 235 | 233 | 231 | 229 | 227 | 225 | 223 |
| 67 | ∞ | ∞ | ∞ | ∞ | ∞ | ∞ | ∞ | 399 | 357 | 317 | 273 | 243 | 241 | 239 | 237 | 235 | 233 | 231 | 229 | 227 |
| 68 | ∞ | ∞ | ∞ | ∞ | ∞ | ∞ | ∞ | 409 | 365 | 325 | 281 | 249 | 245 | 243 | 241 | 239 | 237 | 235 | 233 | 231 |
| 69 | ∞ | ∞ | ∞ | ∞ | ∞ | ∞ | ∞ | 419 | 373 | 333 | 289 | 255 | 249 | 247 | 245 | 243 | 241 | 239 | 237 | 235 |
| 70 | ∞ | ∞ | ∞ | ∞ | ∞ | ∞ | ∞ | 429 | 381 | 341 | 297 | 261 | 253 | 251 | 249 | 247 | 245 | 243 | 241 | 239 |
| 71 | ∞ | ∞ | ∞ | ∞ | ∞ | ∞ | ∞ | 439 | 389 | 349 | 305 | 267 | 257 | 255 | 253 | 251 | 249 | 247 | 245 | 243 |
| 72 | ∞ | ∞ | ∞ | ∞ | ∞ | ∞ | ∞ | 449 | 397 | 357 | 313 | 273 | 261 | 259 | 257 | 255 | 253 | 251 | 249 | 247 |
| 73 | ∞ | ∞ | ∞ | ∞ | ∞ | ∞ | ∞ | 459 | 405 | 365 | 321 | 279 | 265 | 263 | 261 | 259 | 257 | 255 | 253 | 251 |
| 74 | ∞ | ∞ | ∞ | ∞ | ∞ | ∞ | ∞ | 469 | 413 | 373 | 329 | 285 | 269 | 267 | 265 | 263 | 261 | 259 | 257 | 255 |
| 75 | ∞ | ∞ | ∞ | ∞ | ∞ | ∞ | ∞ | 481 | 421 | 381 | 337 | 291 | 273 | 271 | 269 | 267 | 265 | 263 | 261 | 259 |
| 76 | ∞ | ∞ | ∞ | ∞ | ∞ | ∞ | ∞ | 493 | 429 | 389 | 345 | 297 | 277 | 275 | 273 | 271 | 269 | 267 | 265 | 263 |
| 77 | ∞ | ∞ | ∞ | ∞ | ∞ | ∞ | ∞ | 505 | 437 | 397 | 353 | 305 | 281 | 279 | 277 | 275 | 273 | 271 | 269 | 267 |
| 78 | ∞ | ∞ | ∞ | ∞ | ∞ | ∞ | ∞ | 517 | 445 | 405 | 361 | 313 | 285 | 283 | 281 | 279 | 277 | 275 | 273 | 271 |
| 79 | ∞ | ∞ | ∞ | ∞ | ∞ | ∞ | ∞ | 529 | 453 | 413 | 369 | 321 | 289 | 287 | 285 | 283 | 281 | 279 | 277 | 275 |
| 80 | ∞ | ∞ | ∞ | ∞ | ∞ | ∞ | ∞ | 541 | 461 | 421 | 377 | 329 | 295 | 291 | 289 | 287 | 285 | 283 | 281 | 279 |
| 81 | ∞ | ∞ | ∞ | ∞ | ∞ | ∞ | ∞ | 553 | 469 | 429 | 385 | 337 | 301 | 295 | 293 | 291 | 289 | 287 | 285 | 283 |
| 82 | ∞ | ∞ | ∞ | ∞ | ∞ | ∞ | ∞ | 565 | 477 | 437 | 393 | 345 | 307 | 299 | 297 | 295 | 293 | 291 | 289 | 287 |
| 83 | ∞ | ∞ | ∞ | ∞ | ∞ | ∞ | ∞ | 579 | 485 | 445 | 401 | 353 | 313 | 303 | 301 | 299 | 297 | 295 | 293 | 291 |
| 84 | ∞ | ∞ | ∞ | ∞ | ∞ | ∞ | ∞ | 593 | 493 | 453 | 409 | 361 | 319 | 307 | 305 | 303 | 301 | 299 | 297 | 295 |
| 85 | ∞ | ∞ | ∞ | ∞ | ∞ | ∞ | ∞ | 607 | 501 | 461 | 417 | 369 | 325 | 311 | 309 | 307 | 305 | 303 | 301 | 299 |
| 86 | ∞ | ∞ | ∞ | ∞ | ∞ | ∞ | ∞ | 621 | 509 | 469 | 425 | 377 | 331 | 315 | 313 | 311 | 309 | 307 | 305 | 303 |
| 87 | ∞ | ∞ | ∞ | ∞ | ∞ | ∞ | ∞ | 635 | 517 | 477 | 433 | 385 | 337 | 319 | 317 | 315 | 313 | 311 | 309 | 307 |
| 88 | ∞ | ∞ | ∞ | ∞ | ∞ | ∞ | ∞ | 649 | 525 | 485 | 441 | 393 | 343 | 323 | 321 | 319 | 317 | 315 | 313 | 311 |
| 89 | ∞ | ∞ | ∞ | ∞ | ∞ | ∞ | ∞ | 663 | 533 | 493 | 449 | 401 | 349 | 327 | 325 | 323 | 321 | 319 | 317 | 315 |
| 90 | ∞ | ∞ | ∞ | ∞ | ∞ | ∞ | ∞ | 677 | 541 | 501 | 457 | 409 | 357 | 331 | 329 | 327 | 325 | 323 | 321 | 319 |
| 91 | ∞ | ∞ | ∞ | ∞ | ∞ | ∞ | ∞ | 691 | 549 | 509 | 465 | 417 | 365 | 335 | 333 | 331 | 329 | 327 | 325 | 323 |
| 92 | ∞ | ∞ | ∞ | ∞ | ∞ | ∞ | ∞ | 705 | 557 | 517 | 473 | 425 | 373 | 339 | 337 | 335 | 333 | 331 | 329 | 327 |
| 93 | ∞ | ∞ | ∞ | ∞ | ∞ | ∞ | ∞ | 721 | 565 | 525 | 481 | 433 | 381 | 345 | 341 | 339 | 337 | 335 | 333 | 331 |
| 94 | ∞ | ∞ | ∞ | ∞ | ∞ | ∞ | ∞ | 737 | 575 | 533 | 489 | 441 | 389 | 351 | 345 | 343 | 341 | 339 | 337 | 335 |
| 95 | ∞ | ∞ | ∞ | ∞ | ∞ | ∞ | ∞ | 753 | 585 | 541 | 497 | 449 | 397 | 357 | 349 | 347 | 345 | 343 | 341 | 339 |
| 96 | ∞ | ∞ | ∞ | ∞ | ∞ | ∞ | ∞ | 769 | 595 | 549 | 505 | 457 | 405 | 363 | 353 | 351 | 349 | 347 | 345 | 343 |
| 97 | ∞ | ∞ | ∞ | ∞ | ∞ | ∞ | ∞ | 785 | 605 | 557 | 513 | 465 | 413 | 369 | 357 | 355 | 353 | 351 | 349 | 347 |
| 98 | ∞ | ∞ | ∞ | ∞ | ∞ | ∞ | ∞ | 801 | 615 | 565 | 521 | 473 | 421 | 375 | 361 | 359 | 357 | 355 | 353 | 351 |
| 99 | ∞ | ∞ | ∞ | ∞ | ∞ | ∞ | ∞ | 817 | 625 | 573 | 529 | 481 | 429 | 381 | 365 | 363 | 361 | 359 | 357 | 355 |
| 100 | ∞ | ∞ | ∞ | ∞ | ∞ | ∞ | ∞ | 833 | 635 | 581 | 537 | 489 | 437 | 387 | 369 | 367 | 365 | 363 | 361 | 359 |

Table 2. Table of $F(n,S)$ for $51 \leq n \leq 100$ and $1 \leq S \leq 20$.

**Proof.** Except for boundary cases, the result is obtained by comparing the expression in the minimum of the recursion of $F(n, S)$ for $m = m(n, S), m = m(n, S) - 1$ and $m = m(n, S) + 1$. The boundary cases are readily checked. ∎

**Lemma 2.5** *The following hold:*

(i) $\Delta(n, S)$ *is nondecreasing in* $n$.

(ii) *Let* $n \geq 1$ *and* $S \geq 2$. *If* $\Delta(n - m(n, S), S - 1) \leq \Delta(m(n, S), S) + \Delta(m(n, S), S - 1)$, *then* $m(n+1, S) = m(n, S)$. *If not, then* $m(n+1, S) = m(n, S) + 1$.

(iii) *Let* $n \geq 1$ *and* $S \geq 2$. *Then either* $\Delta(n, S) = \Delta(n - m(n, S), S - 1)$ *or* $\Delta(n, S) = \Delta(m(n, S), S) + \Delta(m(n, S), S - 1)$, *whichever is smaller.*

**Proof.** We prove the lemma by induction on $S$ and $n$. For $S = 1$, $\Delta(n, S)$ is nondecreasing in $n$ by inspection. Let $S > 1$. We have $m(0, S) = m(1, S) = 0$ by definition. Note that $m(2, S) = 1$. $\Delta(0, S) = 0$ and $\Delta(1, S) = 2$ which gives the statement for $n = 0$ and $n = 1$. Let $n \geq 2$. By induction, $\Delta(n', S')$ is nondecreasing for $S' = S - 1$ and for $S' = S$ and $n' < n$. Let $m = m(n, S)$. Using Lemma 2.4 we get for $m' < m$

$$\begin{aligned}
\Delta(n + 1 - m' - 1, S - 1) &\geq \Delta(n - m, S - 1) \\
&> \Delta(m - 1, S) + \Delta(m - 1, S - 1) \\
&\geq \Delta(m', S) + \Delta(m', S - 1),
\end{aligned}$$

which means that $m' \neq m(n + 1, S)$. Similarly, for $m' > m + 1$,

$$\begin{aligned}
\Delta(n + 1 - m', S - 1) &\leq \Delta(n - m - 1, S - 1) \\
&\leq \Delta(m, S) + \Delta(m, S - 1) \\
&\leq \Delta(m' - 1, S) + \Delta(m' - 1, S - 1),
\end{aligned}$$

so that $m' \neq m(n + 1, S)$. This leaves only $m(n + 1, S) = m$ or $m(n + 1, S) = m + 1$. The former holds iff $\Delta(n - m, S - 1) \leq \Delta(m, S) + \Delta(m, S - 1)$. In this case, $\Delta(n, S) = \Delta(n - m, S - 1)$ and otherwise $\Delta(n, S) = \Delta(m, S) + \Delta(m, S - 1)$. In the former case, by induction $\Delta(n - 1, S) = \Delta(n - 1 - m, S - 1)$ or $\Delta(n - 1, S) = \Delta(m - 1, S) + \Delta(m - 1, S - 1)$, whichever is smaller, so that by monotonicity, $\Delta(n, S) \geq \Delta(n - 1, S)$. In the latter case, $\Delta(n - 1, S) = \Delta(n - 1 - (m - 1), S - 1)$ or $\Delta(n-1, S) = \Delta(m-2, S) + \Delta(m-2, S-1)$, whichever is smaller, and again it is clear that $\Delta(n, S) \geq \Delta(n-1, S)$. ∎

**Lemma 2.6** *Let $m = m(n, S)$. Then*

$$\Delta(n - m - 1, S - 1) \leq \Delta(n, S) \leq \Delta(n - m, S - 1),$$
$$\Delta(m - 1, S) + \Delta(m - 1, S - 1) \leq \Delta(n, S) \leq \Delta(m, S) + \Delta(m, S - 1).$$

**Proof.** This is an immediate consequence of Lemmas 2.4 and 2.5. ∎

**Lemma 2.7** $\Delta(n, S)$ *is nonincreasing in $S$.*

**Proof.** The lemma follows from (i) of Lemma 2.5, Lemma 2.6 and by induction. Let $m = m(n, S+1)$ and $m' = m(n, S)$. We have

$$\begin{aligned} \Delta(n, S+1) &= \min(\Delta(n-m, S+1), \Delta(m, S+1) + \Delta(m, S)) \\ &\leq \min(\Delta(n-m, S), \Delta(m, S) + \Delta(m, S-1)). \end{aligned}$$

If $m' = m$, then the last expression is $\Delta(n, S)$ and we are done. If $m' > m$, then by Lemma 2.6,

$$\begin{aligned} \Delta(n, S) &\geq \Delta(m'-1, S) + \Delta(m'-1, S-1) \\ &\geq \Delta(m, S) + \Delta(m, S-1). \end{aligned}$$

Similarly, for $m' < m$,

$$\begin{aligned} \Delta(n, S) &\geq \Delta(n - m' - 1, S - 1) \\ &\geq \Delta(n - m, S). \end{aligned}$$

In either case, the desired inequality follows. ∎

For $k \geq 0$, let $x_{k,S} = \min(n \mid \Delta(n, S) > 2^k)$. This is well defined and nondecreasing in $k$. We have $x_{0,S} = 1$, $x_{1,S} = S$, and for $k \geq 0$, $x_{k,1} = 1$.

**Lemma 2.8** *For $S \geq 2$ and $k \geq 0$ we have*

$$\begin{aligned} x_{k+1,S} &\geq x_{k+1,S-1} + x_{k,S-1}, \\ x_{k+1,S} &\leq x_{k+1,S-1} + x_{k,S}, \\ x_{k+1,S} &\leq 2 x_{k+1,S-1}. \end{aligned}$$

**Proof.** The proof is by induction on $S$ and $k$. By Lemma 2.6, if either $\Delta(n - m, S - 1) \leq 2^{k+1}$ or $\Delta(m, S - 1) \leq 2^k$, then $\Delta(n, S) \leq 2^{k+1}$. This must occur if $n < x_{k+1,S-1} + x_{k,S-1}$ so that $x_{k+1,S} \geq x_{k+1,S-1} + x_{k,S-1}$.

Suppose that $x_{k,S} > x_{k+1,S-1}$ and $n \geq 2x_{k+1,S-1}$. Since $m = m(n, S) \leq \lfloor n/2 \rfloor$, $n - m \geq x_{k+1,S-1}$. If $m = n - m = x_{k+1,S-1}$, then $\Delta(n - m, S - 1) > 2^{k+1}$ and $\Delta(m, S) + \Delta(m, S - 1) \geq \Delta(m, S - 1) > 2^{k+1}$ which gives $\Delta(n, S) > 2^{k+1}$. If $m < x_{k+1,S-1}$, then $n - m > x_{k+1,S-1}$ and by Lemma 2.6, $2^{k+1} < \Delta(n - m - 1, S - 1) \leq \Delta(n, S)$, so the last inequality follows.

Suppose that $x_{k,S} \leq x_{k+1,S-1}$ and $n \geq x_{k,S} + x_{k+1,S-1}$. If $m = m(n, S) > x_{k,S}$, then the second inequality of Lemma 2.6 and monotonicity of $\Delta$ imply $\Delta(n, S) > 2^{k+1}$. If $m < x_{k,S}$, then $n - m > x_{k+1,S-1}$ and the first inequality of Lemma 2.6 can be applied to show that $\Delta(n, S) > 2^{k+1}$. Suppose that $m = x_{k,S}$ and $n = x_{k+1,S-1}$. Then both $\Delta(n - m, S - 1) > 2^{k+1}$ and $\Delta(m, S - 1) \geq \Delta(m, S) > 2^k$ and again $\Delta(n, S) > 2^{k+1}$. This completes the proof of the lemma. ∎

**Lemma 2.9** *For $k \geq 0$ and $S \geq 2$, let*

$$x_{k,S}^{(l)} = \sum_{i=0}^{k} \binom{S-1}{i},$$

$$x_{k,S}^{(u)} = \min\left(\binom{S+k-1}{k}, 2^{S-1}\right).$$

*Then $x_{k,S}^{(l)} \leq x_{k,S} \leq x_{k,S}^{(u)}$.*

**Proof.** It is straightforward to check that $x^{(l)}$ satisfies the recursion $x_{k+1,S}^{(l)} = x_{k+1,S-1}^{(l)} + x_{k,S-1}^{(l)}$ with the initial condition $x_{1,S}^{(l)} = 1$. The function $y_{k,S} = \binom{S+k-1}{k}$ satisfies the recursion $y_{k+1,S} = y_{k+1,S-1} + y_{k,S}$ with the correct initial condition. The bound $2^{S-1}$ in the definition of $x^{(u)}$ comes from Theorem 2.3. The result follows by Lemma 2.8. ∎

Lemma 2.9 cannot be used to obtain "nice" expressions for $F(n, S)$ but yields lower and upper bounds by summing $\Delta(n, S)$.

**Lemma 2.10** *For $1 \leq k \leq S - 1$ and $S \geq 2$*

$$F(x_{k,S}^{(l)}, S) \leq \sum_{i=0}^{k} \binom{S-1}{i} 2^{i+1},$$

$$F(x_{k,S}^{(u)}, S) \geq \sum_{i=0}^{k} \binom{S+i-2}{i}(2^i + 1).$$

**Proof.** The inequalities are obtained by summing $\Delta(n, S)$ and making use of Lemma 2.9, ignoring the second case for $x_{k,S}^{(u)}$. Observe that $\binom{S+k-1}{k} = \sum_{i=0}^{k} \binom{S+i-2}{i}$. ∎

## 2.4 Asymptotic bounds

For asymptotics where $S = O(\log(n))$, let

$$f(\gamma, S) = \frac{1}{S} \log F\left(\lfloor 2^{\gamma S} \rfloor, S\right).$$

All logarithms are base 2. The floor/ceiling notation will be omitted (or absorbed into the 'o' and 'O' notation) for real numbers unless needed for clarity. The definition of $f(\gamma, S)$ implies that $F(n, \log(n)/\gamma) = n^{f(\gamma, \log(n)/\gamma)/\gamma + o(1)}$.

For $0 \leq \gamma \leq 1$, let $H(\gamma)$ denote the usual entropy function

$$H(\gamma) = -\gamma \log(\gamma) - (1 - \gamma) \log(1 - \gamma).$$

**Theorem 2.11** *Let $\gamma \leq 1/2$. Then*

$$f(H(\gamma), S) \leq \gamma + H(\gamma) + o(1),$$
$$f((1+\gamma)H(\gamma/(1+\gamma)), S) \geq \gamma + (1+\gamma)H(\gamma/(1+\gamma)) + o(1),$$

*where the $o(1)$ are with respect to $S \to \infty$.*

The proof of the theorem is below. First we make a few observations. For small $\gamma$, the two expressions in the bounds are the same to within $O(\gamma^2)$.

$$\begin{aligned}
H(\gamma) &= -\gamma \log(\gamma) + \log(e)\gamma + O(\gamma^2) \\
&= \gamma(\log(e/\gamma)) + O(\gamma^2), \\
(1+\gamma)H(\gamma/(1+\gamma)) &= -\gamma \log(\gamma) + (1+\gamma) \log(1+\gamma) \\
&= \gamma(\log(e/\gamma)) + O(\gamma^2).
\end{aligned}$$

So for small $\gamma$, the bounds are asymptotically tight. The upper bound is an asymptotic version of that given in [2] with a slightly better constant. More specifically, let $\epsilon = \frac{1}{\log(e/\gamma)}$. Then $\gamma = e2^{-1/\epsilon}$ and $H(\gamma) = (1+\gamma)H(\gamma/(1+\gamma)) + O(\gamma^2) = \frac{\gamma}{\epsilon} + O(\gamma^2)$. Substituting this in the theorem gives

$$F\left(n, \frac{\epsilon}{e}2^{1/\epsilon}\log(n)\right) = n^{1+\epsilon+2^{-O_\epsilon(1/\epsilon)}+o_n(1)},$$

where we used suffixes on the '$O$' and '$o$' notation to clarify the limit involved.

**Proof of Theorem 2.11.** A consequence of the Stirling approximation of $n!$ is that for $\alpha$ constant, $\alpha \leq 1/2$,

$$\sum_{i=0}^{\alpha S} \binom{S}{i} = O\left(\binom{S}{\alpha S}\right) = 2^{H(\alpha)S(1+o(1))}.$$

Thus, for $n = 2^{H(\gamma)S}$, $\Delta(n, S) \leq 2^{\gamma S(1+o(1))}$. Using these observations and Lemma 2.10, one can estimate

$$\begin{aligned} F(n, S) &\leq \binom{S}{\gamma S} 2^{\gamma S(1+o(1))} \\ &\leq 2^{(\gamma+H(\gamma))S(1+o(1))}. \end{aligned}$$

Similarly, for $\alpha$ constant,

$$\sum_{i=0}^{\alpha S} \binom{S+i}{i} = O\left(\binom{S(1+\alpha)}{\alpha S}\right) = 2^{(1+\alpha)H(\alpha/(1+\alpha))S(1+o(1))},$$

so that for $n = 2^{(1+\gamma)H(\gamma/(1+\gamma))}$, $\Delta(n, S) \geq 2^{\gamma S(1+o(1))}$. This gives the second inequality. ∎

## 2.5 Asymptotics of $TS$.

The minimum time-space product is asymptotically determined by the next theorem.

**Theorem 2.12**

$$TS(n) = 2^{2\sqrt{\log(n)}(1+o(1))}n.$$

**Proof.** To show that $TS(n) \leq 2^{2\sqrt{\log(n)}(1+o(1))}$, let $S = 2^{\sqrt{\log(n)}}$. Then

$$\begin{aligned} x_{c\sqrt{\log(n)},S} &\geq \sum_{i=0}^{c\sqrt{\log(n)}} \binom{S}{i} \\ &\geq \binom{2^{\sqrt{\log(n)}}}{c\sqrt{\log(n)}} \\ &\geq n^{c+o(1)}. \end{aligned}$$

This implies that $\Delta(n,S) \leq 2^{\sqrt{\log(n)}(1+o(1))}$ and

$$F(n,S) \leq 2^{\sqrt{\log(n)}(1+o(1))}n.$$

For the lower bound, we estimate $F(n,S)$ for $S = 2^{d\sqrt{\log(n)}}$ from below.

$$\begin{aligned} x_{c\sqrt{\log(n)},S} &\leq \binom{S + c\sqrt{\log(n)} - 1}{c\sqrt{\log(n)}} \\ &\leq \left(S + c\sqrt{\log(n)}\right)^{c\sqrt{\log(n)}} \\ &\leq 2^{cd\log(n)(1+o(1))} = n^{cd(1+o(1))}. \end{aligned}$$

This implies that $\Delta(n,S) \geq 2^{\frac{\sqrt{\log(n)}}{d}(1+o(1))}$ and

$$\begin{aligned} F(n,S)S &\geq 2^{\sqrt{\log(n)}(d+\frac{1}{d})(1+o(1))}n(1-o(1)) \\ &\geq 2^{2\sqrt{\log(n)}(1+o(1))}n, \end{aligned}$$

since $d + \frac{1}{d} \geq 2$ for all $d > 0$. We used the observation that the last term of the sum $\sum_{i=0}^{k} \binom{S+i-2}{i}$ asymptotically dominates it. ∎

## A  Appendix

### A.1  Mathematica expressions for dynamic evaluation of $F(n,S)$

```
FBen[1,S_] := If[S >= 1, 1, Infinity];
FBen[n_,0] := Infinity;

FBen[n_,S_] := FBen[n,S] = Module[{m},
    Min[Table[FBen[m,S]+FBen[m,S-1]+FBen[n-m,S-1],{m,1,n-1}]]
  ];
```

If the goal is to compute a table, then Lemma 2.5 can be used for efficient computation by keeping track of $m(n,S)$.

```
FBenTable[nmax_,Smax_] := Module[{m,n,nt,S,T,Tt,t1,t2,t},
    T = {Join[{1}, Table[Infinity,{n,2,nmax}]]};
    For[S=2, S<=Smax, S++,
      m=1;
```

```
      Tt = {1,3};
      For[n=3, n<=nmax, n++,
        t1=Tt[[m]]+T[[S-1,m]]+T[[S-1,n-m]];
        t2=Tt[[m+1]]+T[[S-1,m+1]]+T[[S-1,n-m-1]];
        If[t2<t1,
          m++;
          t = t2;
          ,
          t = t1;
        ];
        Tt = Append[Tt,t];
        If[t==Infinity,
          Tt = Join[Tt,Table[Infinity,{nt,n+1,nmax}]];
          Break[];
        ,Null]
      ];
      T = Append[T,Tt];
    ];
    Return[T];
  ];
```